\documentclass[12pt,reqno]{article}

\usepackage[top=2.5cm,bottom=2.5cm,left=2.5cm,right=2.5cm]{geometry}
\usepackage{graphicx}
\usepackage{latexsym}
\usepackage{amsmath,amssymb}
\usepackage{amscd}
\usepackage[arrow,matrix]{xy}
\usepackage{graphicx}
\usepackage{xcolor}
\usepackage{times}
\usepackage[]{subfigure}
\RequirePackage{hyperref}
\RequirePackage{bookmark}
 \usepackage{tabu}
 \usepackage{cite}
 \usepackage{float}

\usepackage{amsfonts}

\usepackage{amscd,amsmath,amsthm,amssymb}

\theoremstyle{plain}
\newtheorem{theorem}{Theorem}[section]

\theoremstyle{definition}
\newtheorem{remark}[theorem]{Remark}

\begin{document}

\title{\textbf{Surfaces in a strict Walker 3-manifold that contain non-null curves with zero torsion
}}
\author {{El Hadji Baye Camara$^{1}$}\thanks{{
 E-mail: \texttt{elhadjibaye.camara@ucad.edu.sn } (E. B. Camara)}},
 \texttt{} Athoumane Niang$^{2}$\footnote{{
 E-mail: \texttt{athoumane.niang@ucad.edu.sn} (A. Niang)}},
 \texttt{} Ameth Ndiaye$^{3}$\footnote{{
 E-mail: \texttt{ameth1.ndiaye@ucad.edu.sn} (A. Ndiaye)}},\texttt{ } Adama Thiandoum$^{4}$\footnote{{
 E-mail: \texttt{adama1.thiandoum@ucad.edu.sn} (A. Thiandoum)}} \\
\begin{small}{$^{1,2,4}$D\'epartement de Mathematiques et Informatique, FST, UCAD, Dakar, S\'en\'egal.}\end{small}\\
\begin{small}{$^{3}$D\'epartement de Math\'ematiques, FASTEF, UCAD, Dakar, Senegal.}\end{small}}
\date{}
\maketitle%
\textbf{}
\begin{abstract} 
Given a non-null curve $\gamma$ in a strict Walker 3-manifold, first we show that (locally) $\gamma$ lies in a flat cylinder with a null axis.
Secondly, we construct an example of such a curve $\gamma$ and such a cylinder $S$ that contains $\gamma$.
In particular, the hypothesis that $S$ is totally geodesic has some consequence on the geometry of the ambient Walker 3-manifold.
\end{abstract}
\begin{small} {\textbf{MSC 2020:}  53A10, 53C42, 53C50.}
\end{small}\\
\begin{small} {\textbf{Keywords:} Curve, Surface, Walker manifold, Curvature, Mean curvature, Gauss curvature.} 
\end{small}\\

\section{Introduction}
The three Walker manifolds are described in terms of suitable coordinates $(x,y,z)$ of the manifold $\mathbb{R}^3$ and their metric depends on an arbitrary function $f = f(x,y,z)$ and are given by:
\begin{eqnarray*}
g_f = \varepsilon\, dy^2 + 2\, dx\, dz + f\, dz^2 \quad, \quad \varepsilon = \pm 1.
\end{eqnarray*}
These manifolds are denoted by $(M, g_f)$ or $M_f$ (see \cite{Walker}).\\
When the function $f$ depends only on the variables $y$ and $z$, $M_f$ is called a strict Walker 3-manifold.\\
The metric tensor $g_f$ is given by the matrix:

\begin{align} \label{1}
A = \begin{pmatrix}
0 & 0 & 1 \\
0 & \varepsilon & 0 \\
1 & 0 & f
\end{pmatrix}.
\end{align}.

Hence, if $(M, g_f)$ is a strict Walker manifold $f(x, y, z) = f(y, z)$, then the associated Levi-Civita connection satisfies that the non-zero components of its Christoffel symbols are given by

\begin{align} \label{2}
\nabla_{\partial_y} \partial_z = \frac{1}{2} f_y\, \partial_x
\quad \text{and} \quad
\nabla_{\partial_z} \partial_z = \frac{1}{2} f_z\, \partial_x - \frac{\varepsilon}{2} f_y\, \partial_y .
\end{align}

And the non-zero components of the curvature of the strict Walker 3-manifold $(M, g_f)$ are given by,

\begin{align} \label{3}
\left\{
\begin{array} {llllllll}
R(\partial_y, \partial_z)\, \partial_y &= \frac{1}{2} f_{yy}\, \partial_x \cr
                                                                           \cr
R(\partial_y, \partial_z)\, \partial_z &= -\frac{\varepsilon}{2} f_{yy}\, \partial_y .
\end{array}
\right.
\end{align}

Denoted by $(e_1, e_2, e_3)$ the canonical basis in $\mathbb{R}^3$ associated to the coordinate $(x, y, z)$.

If $u = (u_1, u_2, u_3)$ and $v = (v_1, v_2, v_3)$ are two vectors in $\mathbb{R}^3$, we define the vector product by

\begin{eqnarray*}
g_f(u \times_{f} v, w) = \det(u, v, w).
\end{eqnarray*}

One finds the following formula:

\begin{equation} \label{4}
u \times_{f} v =
\Big(\begin{vmatrix}
u_1 & v_1 \\
u_2 & v_2 
\end{vmatrix}
-f \begin{vmatrix}
u_2 & v_2 \\
u_3 & v_3 
\end{vmatrix}\Big) e_1
- \varepsilon
\begin{vmatrix}
u_1 & v_1 \\
u_3 & v_3
\end{vmatrix} e_2
+
\begin{vmatrix}
u_2 & v_2 \\
u_3 & v_3
\end{vmatrix} e_3,
\quad \text{(see \cite{Camara, Niang})}.
\end{equation}

We will denote $g_f$ by $\langle \cdot, \cdot \rangle$.

Let us introduce an orthonormal frame $\left(T(s), N(s), B(s)\right)$ along a non-null curve $\gamma(s)$ parametrized by arc length $s$. The first unit tangent vector field $T(s)$ is defined by $T(s) = \dot{\gamma}(s)$.

The acceleration vector field $\nabla_{\dot{\gamma}} T$ is non-null, that is non-lightlike. Denote by $\varepsilon_2 = \pm 1$ the sign of $\langle \nabla_{\dot{\gamma}} T,\; \nabla_{\dot{\gamma}} T \rangle$.

The curvature $\kappa(s)$ of $\gamma$ is
\begin{eqnarray} \label{5}
\kappa(s) = \sqrt{\varepsilon_2 \langle \nabla_T T,\; \nabla_T T \rangle} > 0.
\end{eqnarray}

Then we can find a unit vector $N(s)$ on $\gamma$ such that
\begin{eqnarray*}
\nabla_T T = \varepsilon_2 \kappa N
\quad \text{and} \quad \langle N, N \rangle = \varepsilon_2.
\end{eqnarray*}

We define a unit vector field $B$ by $B = \varepsilon_3 T \times_{f} N$, where $\varepsilon_3 = \langle B, B \rangle$.

Thus we obtains the following Frenet-Serret equations:

\begin{align} \label{6}
\left\{
\begin{array} {lllll}
\nabla_T T &= \varepsilon_2 \kappa N \\
                                      \\
\nabla_T N &= -\varepsilon_1 \kappa T - \varepsilon_3 \tau B \\
                                      \\
\nabla_T B &= \varepsilon_2 \tau N .
\end{array}
\right.
\end{align}

The function $\tau$ is called the torsion of $\gamma$, see \cite{Erjavec}.\\
By the Remark 1 in \cite{Erjavec}, when the ambient manifold $M$ is of constant curvature, $\tau = 0$ is equivalent to the condition that $\gamma$ lies in a totally geodesic surface of $M$.\\ In \cite{Carmo} and \cite{Lopez}, it is shown that in the Euclidean 3-space or in the Minkowski 3-space, $\tau = 0$ is equivalent to the condition that $\gamma$ lies in an affine plane.

\vspace{0.5cm}

The aim of this work is to study the analogue of this result in a strict Walker 3-manifold $M_f$.\\
Our main result is the following:

\begin{theorem} \label{T1}
Let $\gamma$ be a non-null curve in $M_f$ with zero torsion. Then $\gamma$ lies (locally) in a flat cylinder with a null axis.
And such a cylinder, being totally geodesic, has some consequence on the geometry of the ambient manifold $M_f$.
More over we give example zero torsion curve lying in a such cylinder.
\end{theorem}
\vspace{0.2cm}
The paper is organised as follows.\\
The section 2 contains notations and some recalls and formulas for curves and surfaces in $M_f$.\\
The section 3 contains the proof of theorem \ref{T1}. We end Section 3 by the study of an example of zero torsion curve $\gamma$ in $M_f$ and a cylinder $S$ that contains $\gamma$.

\section{Some recalls and formulas for curves and surface}
Let $(M, g_{f})$ be a pseudo-Riemannian manifold with its canonical Levi-Civita connection denoted by $\nabla$. Assume that $M$ has dimension 3, and let $\Sigma \subset M$ be a semi-Riemannian surface of $(M, g_f)$.\\
Denote by $\eta$ a unit normal vector field on $\Sigma$ with sign $\delta = \pm 1$.\\
The Levi-Civita connection $D$ of the induced metric on $\Sigma \subset M$, and $\nabla$ are related by:\\
if $X$, $Y$ and $Z$ are tangent vector fields on $\Sigma,$ we have:

\begin{eqnarray} \label{7}
\nabla_X Y = D_X Y + h(X,Y)\eta .
\end{eqnarray}
where $h$ is the second fundamental form.\\

The second fundamental form $h$ and Weingarten map $S$ of $\Sigma$ are related by:
\begin{align} \label{8}
\left\lbrace
\begin{array}{llllll}
-\nabla_X \eta &= SX \cr
g_{f}(SX, Y) &= \delta h(X,Y) = g_{f}(\nabla_X Y, \eta).
\end{array}
\right.
\end{align}

If we denote by $R^M$ the curvature of $(M, g_{f})$ and $R$ the curvature of $(\Sigma, i^*\gamma)$ for $i: \Sigma \hookrightarrow M$ the inclusion map, we have

\begin{align} \label{9}
\langle R(X,Y)Z, W \rangle &= \langle R^M(X,Y)Z, W \rangle +\delta\big( h(X,W)h(Y,Z) - h(Y,W)h(X,Z)\big).
\end{align}
for $X, Y, Z, W$ tangent to $\Sigma$.\\

And the sectional curvature $K$ on $\Sigma$ and $K^M$ of $M$ are related by:
\begin{align} \label{10}
K(X,Y) = K^M(X,Y) + \delta \det S .
\end{align}

Now let $\gamma = \gamma(t)$ be a non-null curve in a strict Walker 3-manifold $(M_f)$.\\
We put $\gamma(t) = \Big(x(t),\; y(t),\; z(t) \Big)$.\\

So
\[
\gamma'(t) = x'(t)\, \partial_x + y'(t)\, \partial_y + z'(t)\, \partial_z,
\]
and
\[
\langle \gamma'(t), \gamma'(t) \rangle = 2x'(t)z'(t) + \varepsilon (y'(t))^2 + \left(z'(t)\right)^2 f(y(t), z(t))
\]

The unit tangent vector $T(t)$ of $\gamma$ is:
\[
T(t) = t_1(t)\, \partial_x + t_2(t)\, \partial_y + t_3(t)\, \partial_z
\]
with
\begin{align} \label{11}
\left\{
\begin{array}{llllllllll}
t_1(t) &= \frac{x'(t)}{\lvert \gamma'(t) \rvert} \\
                                                                                                     \\
t_2(t) &= \frac{y'(t)}{\lvert \gamma'(t) \vert} \\
                                                                                   \\
t_3(t) &= \frac{z'(t)}{\lvert \gamma'(t) \rvert}
\end{array}
\right. 
\end{align} 
where $\lvert \gamma'(t) \rvert = \sqrt{\lvert 2x'(t) z'(t) +\varepsilon (y'(t))^2 + (z'(t))^2 f\big(y(t), z(t)\big)\rvert}$\\

Let $\xi(t) = X(t) \partial x + Y(t) \partial y + Z(t) \partial z$, be a vector field on $\gamma$. \\
By using \eqref{2} an easy computation gives that
\begin{align} \label{12}
\nabla_T \xi = \Bigl( X' + \frac{1}{2} f_y (Y t_3 + Z t_2) + \frac{1}{2} Z t_3 f_z \Bigr)\partial x
+ \Bigl( Y' - \frac{\varepsilon}{2} Z t_3 f_y \Bigr)\partial y + Z' \partial z.
\end{align}

In particular, using column notation for vector $\xi(t) = T(t)$ the unit tangent of $\gamma$, we have 
\begin{align} \label{13}
\nabla_T T =
\begin{pmatrix}
t_1' + f_y t_2 t_3 + \frac{1}{2} t_3^2 f_z\\[6pt]
t_2' - \frac{\varepsilon}{2} t_3^2 f_y\\[6pt]
t_3'
\end{pmatrix}
:=
\begin{pmatrix}
A_1\\[6pt]
A_2\\[6pt]
A_3
\end{pmatrix}. 
\end{align}

That is $\nabla_T T = A_1 \partial x + A_2 \partial y + A_3 \partial z$,
where $A_1, A_2$ and $A_3$ are defined by \eqref{13}.

By \eqref{5}, the curvature $\kappa$ of $\gamma$ satisfies
\begin{align} \label{14}
\kappa^2 = \bigl| 2 A_1 A_3 + \varepsilon A_2^2 + f A_3^2 \bigr| > 0, 
\end{align}
and
\begin{align} \label{15}
\varepsilon_2 B = \frac{1}{\kappa} T \times_f \nabla_T T.
\end{align} (by \eqref{6}).
\bigskip

\noindent
\section{Proof of the main results}
\subsection{Proof of Theorem \ref{T1}}
We keep the notation in Section 2.
We will identify the ambient space \( M_{f} = \mathbb{R}^3 \) with its tangent space \( T_p \mathbb{R}^3 \) of any \( p \in M_{f} \).

We associated to the curve,
\[
\gamma(t) = \big(x(t),\; y(t),\; z(t)\big).
\]
The vector field on \( \gamma \) which is the duplication of \( \gamma(t) \) at the point \( \gamma(t) \), also denoted by \( \gamma(t) \) :
\begin{equation*}
\gamma(t) = x(t) \, \partial_x|_{\gamma(t)} + y(t) \, \partial_y|_{\gamma(t)} + z(t) \, \partial_z|_{\gamma(t)}.
\end{equation*}

We recall that \( A \) denotes the matrix \( (g_{ij}) \) in \eqref{1}.

Let \( \bullet \) denote the Euclidean dot product in the coordinates \( (x, y, z) \).

In particular, we have the identity:
\begin{align} \label{16}
\langle \gamma,\; B \rangle = \gamma \bullet A B .
\end{align}

Let us take the derivative of \eqref{16} with respect to \( T \) and obtain:
\begin{align} \label{17}
\langle \nabla_{T} \gamma,\; B \rangle + \langle \gamma,\; \nabla_T B \rangle = (\gamma \bullet A B)' .
\end{align}

The torsion \( \tau = 0 \) is equivalent to \( \nabla_T B = 0 \).\\
Thus \eqref{17} becomes:
\begin{align*}
\langle \nabla_{T} \gamma,\; B \rangle &= (\gamma \bullet A B)'\cr
&= (\gamma' \bullet A B) + (\gamma \bullet (A B)') \cr
&= \langle \gamma',\; B \rangle + \gamma \bullet (A B)'\cr
&= \langle \lvert\gamma'\rvert T,\; B \rangle + \gamma \bullet (A B)'.
\end{align*}

Since \(\langle T,\; B \rangle = 0\), the relation in \eqref{17} reduces to
\begin{align} \label{18}
\langle \nabla_T \gamma,\; B \rangle = (\gamma \bullet AB)'.
\end{align}

But
\begin{align*}
\begin{array} {llllll}
\nabla_{T} \gamma &= \gamma' + x(t)\nabla_T \partial_x + y(t)\nabla_T \partial_y + z(t)\nabla_T \partial_z \cr
&= |\gamma'| T + y(t)\nabla_T \partial_y + z(t)\nabla_T \partial_z .
\end{array}
\end{align*}

Since \(\nabla_T \partial_x = 0 \quad \text{(by (2))}\).

Therefore, when the torsion \(\tau = 0\), \eqref{17} gives us:
\begin{align} \label{19}
y(t) \langle \nabla_T \partial_y,\; B \rangle + z(t) \langle \nabla_T \partial_y,\; B \rangle = \gamma \bullet (AB) .'
\end{align}

Since \(\gamma = \gamma(t)\) is a non-null curve, then
\[
\langle \gamma'(t),\; \gamma'(t) \rangle = 2x'(t)z'(t) + \varepsilon (y'(t))^2 + (z'(t))^2 f\big(y(t), z(t)\big).
\]
is not zero at any point \(t_0\).

Thus \(y'(t_0) \neq 0\) or \(z'(t_0) \neq 0\).

This shows that locally, the curve can be parametrized by the coordinate \(y\) or \(z\).

Hence, one can parametrize \(\gamma\) locally as follows:

\begin{itemize}
\item[(a)] \(\gamma : y \mapsto \Big(x(y),\; y,\; z(y)\Big)\) ;
\item[(b)] \(\gamma : z \mapsto \Big(x(z),\; y(z),\; z \Big)\).
\end{itemize}

Assume first that \(\gamma\) is given by
\[
\gamma(y) = \Big(x(y),\; y,\; z(y)\Big).
\]

In view of \eqref{19}, we define the three functions \(C_1,\; C_2,\; C_3\) of the variable \(y\) by

\begin{align} \label{20}
\left\{
\begin{array}{lllllll}
C_1(y) &= \langle \nabla_T \partial_y,\; B(y) \rangle \cr
                                                   \cr
C_2(y) &= \langle \nabla_T \partial_z,\; B(y) \rangle \cr
                                                   \cr
C_3(y) &= \gamma(y) \bullet (AB)'(y).
\end{array}
\right. 
\end{align}

From \eqref{19}, we see that \(\gamma\) lies on the surface
\begin{eqnarray} \label{21}
y \, C_1(y) + z \, C_2(y) = C_3(y).
\end{eqnarray}

The equation \eqref{21} defines a cylinder \(S\) whose axis is colinear with \(\partial x\). Clearly, \(S\) contains \(\gamma\).

If \(\gamma\) has the form (b), then equation \eqref{19} takes the form
\begin{eqnarray} \label{22}
y \, C_1(z) + z \, C_2(z) = C_3(z) .
\end{eqnarray}

which is a cylinder with null axis colinear with \(\partial x\).

In general, the surface defined by \eqref{21} or \eqref{22} has the form
\[
F(y,z) = 0
\]
and locally takes one of the following two forms:

\textbf{Case 1:} \quad \(S : (x,y) \xrightarrow{\Phi} x(1,0,0) + \Big(0, y, \varphi(y)\Big)\)

\textbf{Case 2:} \quad \(S : (x,z) \xrightarrow{\Phi} x(1,0,0) + \Big(0, \varphi(z), z\Big)\)

Let us show that, in each \textbf{Case 1} or \textbf{Case 2}, the surface \(S\) is flat.

We consider the \textbf{Case 1}:

\begin{equation*}
\Phi(x,y) = x\big(1, 0, 0\big) + \Big(0, y, \varphi(y)\Big).
\end{equation*}

Thus, we have
\[
\Phi_x = \partial x \quad \text{and} \quad \Phi_y = \partial y + \varphi'(y)\, \partial z .
\]

And by \eqref{2}, we get

\begin{align} \label{23}
\left\{
\begin{array}{lllllllll}
\Phi_{xx} &= 0 \cr
\Phi_{xy} &= 0 \cr
\Phi_{yy} &= \varphi'(y)\Big(f_{y} + \frac{1}{2} \varphi'(y) f_z \Big)\partial_x -\frac{\varepsilon}{2}\varphi'^{2}(y) f_y\partial_y + \varphi''(y) \partial_z .
\end{array}
\right.
\end{align}

By \eqref{23} and \eqref{8}, one finds that the determinant of the shape operator is 0. This shows that the surface is extrinsically flat.

And \eqref{2} and the fact that $\Phi_x = \partial x$, we see that the surface is intrinsically flat (by \eqref{10} and \eqref{23}).

Thus $S$ is flat.\\
Similarly, if we consider the \textbf{Case 2}: 

\begin{equation*}
\Phi(x,z) = x\big(1, 0, 0\big) + \Big(0, \varphi(z), z\Big).
\end{equation*}

Thus, we have
\[
\Phi_x = \partial x \quad \text{and} \quad \Phi_z =  \varphi'(y)\, \partial y + \partial z .
\]

And by \eqref{2}, we get

\begin{align} \label{023}
\left\{
\begin{array}{lllllllll}
\Phi_{xx} &= 0 \cr
\Phi_{xz} &= 0 \cr
\Phi_{zz} &= \Big(\varphi'(z) f_{y} + \frac{1}{2} f_z \Big)\partial_x  + \Big(\varphi''(z)  - \frac{\varepsilon}{2} f_y \Big)\partial_y .
\end{array}
\right.
\end{align}

We get easily that the surface $S$ given in the \textbf{Case 2} is flat.\\
The unit normal vector field $\eta$ is given by :

In the Case 1 we have
\begin{align*}
\eta = \frac{\varPhi_x \times_f \varPhi_y}{\lvert\varPhi_x \times_f \varPhi_y \rvert} = \frac{1}{\lvert \varphi'(y) \rvert} \Big( 1,\; -\varepsilon\varphi'(y),\; 0 \Big)  \qquad \text{ (by \eqref{4})}.
\end{align*}

Also in the Case 2 we have
\begin{align*}
\eta = \frac{\varPhi_x \times_f \varPhi_z}{\lvert\varPhi_x \times_f \varPhi_z \rvert} =  \Big( \varphi'(z),\; -\varepsilon,\; 0 \Big) \qquad \text{ (by \eqref{4})} .
\end{align*}

Hence,\\
In the \textbf{Case 1}, $S$ is totally geodesic if and anly if,
\begin{equation} \label{00}
\varphi''(y) +\frac{\varepsilon}{2}\varphi'^3(y)f_y\big(y, \varphi(y)\big) = 0 .
\end{equation}
Also in the \textbf{Case 2} an esaly computation shows $S$ is totally geodesic if and anly if,
\begin{equation} \label{01}
\frac{\varepsilon}{2} f_y \big(\varphi(z), z \big)- \varphi''(z) = 0 .
\end{equation}
We have the fillowing remark.\\

\begin{remark}
The equation \eqref{19} in $\mathbb{R}^3$ or in $\mathbb{R}^3_1$ takes the form
\begin{align*}
0 = \gamma \bullet (JB)',
\end{align*}
when
\begin{align*}
J =
\begin{pmatrix}
1 & 0 & 0 \\
0 & 1 & 0 \\
0 & 0 & 1
\end{pmatrix} \quad \text{in } \mathbb{R}^3, \quad \text{and} \quad
J =
\begin{pmatrix}
1 & 0 & 0 \\
0 & 1 & 0 \\
0 & 0 & -1
\end{pmatrix} \quad \text{in } \mathbb{R}^3_1.
\end{align*}
\end{remark}

In these cases $B = B_0$ is a constant vector, and we get:
\[
\gamma \bullet B_0 = \text{const}.
\]

\subsection{Examples of curve with zero torsion in $M_f$}
To find surfaces satisfying \eqref{21} or \eqref{22}, one must find a non-null curve with zero torsion.

The binormal vector $B(t)$ of $\gamma$ satisfies
\[
\varepsilon_2 B(t) = \frac{1}{\kappa}\, T \times_f \nabla_T T
=
\begin{pmatrix}
B_1(t) \\
B_2(t) \\
B_3(t)
\end{pmatrix}.
\]

An easy computation using \eqref{2} gives:
\begin{align} \label{24}
\left\{
\begin{array} {llllll}
\langle \nabla_T \partial_y, B \rangle &= \frac{1}{2} t_3 f_y B_3 \\
                                                                     \\
\langle \nabla_T \partial_z, B \rangle &= -\frac{1}{2} t_3 f_y B_2 +
\left(\frac{1}{2} \partial_z f_y + t_3 B_3 \right) B_3 .
\end{array}
\right.
\end{align}

Since $\langle B, B \rangle \neq 0$, then $B_2(t) \neq 0$ or $B_3(t) \neq 0$ for each $t$ (by \eqref{1}).

Now, $\tau = 0$ if and only if 
\begin{equation*} 
\left( \frac{1}{\kappa} \right)'
\left( T \times_f \nabla_T T \right) + \frac{1}{\kappa} \nabla_T \Big( T \times_f \nabla_T T \Big) = 0 .
\end{equation*} 
We will assume that the curvature $\kappa$ is a positive constant. In this case, curve $\gamma$ has zero torsion if and only if
\begin{equation} \label{25}
\nabla_T \Big( T \times_f \nabla_T T \Big) = 0 .
\end{equation}

We have denoted
\[
T = t_1 \partial x + t_2 \partial y + t_3 \partial z .
\]
and
\[
\nabla_T T = A_1 \partial x + A_2 \partial y + A_3 \partial z
\quad \text{(as in \eqref{11} and \eqref{13})}.
\]

Thus, using \eqref{4}, with \eqref{13}, we get:
\begin{align} \label{26}
T \times_f \nabla_T T =
\begin{pmatrix}
\left| \begin{array}{cc}
t_1 & A_1 \\
t_2 & A_2
\end{array} \right| - f \left| \begin{array}{cc}
t_2 & A_2 \\
t_3 & A_3
\end{array} \right| \\
                                       \\
-\varepsilon \left| \begin{array}{cc}
t_1 & A_1 \\
t_3 & A_3
\end{array} \right| \\
                                       \\
\left| \begin{array}{cc}
t_2 & A_2 \\
t_3 & A_3
\end{array} \right|
\end{pmatrix}
\end{align}

By using the formula in \eqref{12} applied to $T \times_f \nabla_T T$, one gets that $\tau = 0$ and $\kappa = \text{cste}$ gives:

\begin{align} \label{27}
\left\lbrace
\begin{array}{llllllllllllll}
(i) \quad \left(
\left| \begin{array}{cc}
t_1 & A_1 \\
t_2 & A_2
\end{array} \right|
- f \left| \begin{array}{cc}
t_2 & A_2 \\
t_3 & A_3
\end{array} \right|
\right)' + \frac{1}{2} f_y \left(
- \varepsilon t_3 \left| \begin{array}{cc}
t_1 & A_1 \\
t_3 & A_3
\end{array} \right|
+ t_2 \left| \begin{array}{cc}
t_2 & A_2 \\
t_3 & A_3
\end{array} \right|
\right)
+ \frac{1}{2} t_3 f_z \left| \begin{array}{cc}
t_2 & A_2 \\
t_3 & A_3
\end{array} \right| = 0 \\
                        \\
(ii) \quad \left(
\left| \begin{array}{cc}
t_1 & A_1 \\
t_3 & A_3
\end{array} \right|
\right)' - \frac{1}{2} t_3 f_y \left| \begin{array}{cc}
t_2 & A_2 \\
t_3 & A_3
\end{array} \right| = 0 \\
                       \\
(iii) \quad \left(
\left| \begin{array}{cc}
t_2 & A_2 \\
t_3 & A_3
\end{array} \right|
\right)' = 0
\end{array}
\right.
\end{align}

Since $T \times_f \nabla_T T$ is non-null, from \eqref{26}, one sees that
\[
\left| \begin{array}{cc}
t_1 & A_1 \\
t_3 & A_3
\end{array} \right| \neq 0
\quad \text{or} \quad
\left| \begin{array}{cc}
t_2 & A_2 \\
t_3 & A_3
\end{array} \right| \neq 0
\quad \text{at any point} t.
\]

The relation $(iii)$ in \eqref{27} gives that
\[
\left| \begin{array}{cc}
t_2 & A_2 \\
t_3 & A_3
\end{array} \right| = C_1 \quad \text{a constant}.
\]

We will assume $C_1 = 0$, so that $(iii)$ becomes

\begin{align} \label{28}
\begin{pmatrix}
A_2 \\
A_3
\end{pmatrix}
= \lambda
\begin{pmatrix}
t_2 \\
t_3
\end{pmatrix},
\end{align}
for some function $\lambda = \lambda(t)$.

The remaining relations $(ii)$ and $(i)$ of \eqref{27} give:
\begin{align} \label{29}
\left\{
\begin{aligned}
\left| \begin{array}{cc}
t_1 & A_1 \\
t_3 & A_3
\end{array} \right| &= C \quad \text{non-zero constant}, \\
\left( \left| \begin{array}{cc}
t_1 & A_1 \\
t_2 & A_2
\end{array} \right| \right)' &= \frac{\varepsilon}{2} f_y C .
\end{aligned}
\right. 
\end{align}

The curvature $\kappa$ satisfies:
\begin{align} \label{30}
\kappa^2 = \left| 2A_1 A_3 + \varepsilon A_2^2 + A_3^2 f \right| > 0 .
\end{align}

Now we put $t_1 = 0$, so $x(t) = x_0$ a constant real.

The relation \eqref{28} and the first relation in \eqref{29} give us:
\begin{align} \label{31}
A_2 = \lambda t_2, \quad A_3 = \lambda t_3, \quad\text{and} \quad A_1 = -\frac{C}{t_3} 
\end{align}

We put $\gamma = \gamma(y) = \Big(x_0,\; y,\; h(y)\Big)$.

So
\begin{equation} \label{32}
\left\{
\begin{array} {lllll}
t_2(y) &= \frac{1}{\sqrt{\lvert\varepsilon + h'(y)^2 \, f\big(y, h(y)\big)\rvert}} \\
                                                                          \\
t_3(y) &= \frac{h'(y)}{\sqrt{\lvert\varepsilon + h'(y)^2 \, f\big(y, h(y)\big)\rvert}} .
\end{array}
\right.
\end{equation}

So we get
\begin{equation*}
t_3 = h'(y) \, t_2
\end{equation*}

The relations in \eqref{32} show that $t_2$ and $t_3$ are determined by the function $h'(y)$.

Now we use the last relation \eqref{29} with $t_3 = h'(y) \, t_2$ to get
\begin{equation} \label{33}
\left[\frac{1}{h'(y)}\right]' = \frac{\varepsilon}{2} \, f\left(y, h(y)\right) .
\end{equation}

In order to get a solution of \eqref{33}, we will assume that the function $f$ depends only on the variable $y$; that is
\begin{equation*}
f = f(y), \quad \text{with } f \neq 0.
\end{equation*}

Thus, we get
\begin{equation} \label{34}
\frac{1}{h'(y)} = \frac{\varepsilon}{2} \int_{y_0}^{y} f(\tau) \, d\tau + r_0 .
\end{equation}
for a real non-zero $r_0$.

It remains to find $\lambda$. We find $\lambda$ a solution of the equation
\begin{equation*}
\kappa^2 = \left| \frac{-2c}{t_3} t_2 \lambda + \lambda^2 \left( \varepsilon t_2^2 + t_3^2 f(y) \right) \right|.
\end{equation*}

The curves $\gamma = \gamma(y)$ determined by $h = h(y)$ such that,
\begin{equation*}
h'(y) = \frac{1}{\frac{\varepsilon}{2} \int_{y_0}^{y} f(\tau) d\tau + r_0}.
\end{equation*}

Next let us find the surfaces. We put
\[
B =
\begin{pmatrix}
B_1(t) \\
B_2(t) \\
B_3(t)
\end{pmatrix}.
\]

We have also,
\begin{align*}
B &= \frac{1}{\varepsilon_2\kappa} T_{\times_f} \nabla_T T
= \frac{1}{\varepsilon_2\kappa}
\begin{pmatrix}
-t_2 A_1 \\
\varepsilon t_3 A_1 \\
0
\end{pmatrix}.
\end{align*}

By using \eqref{31} with $t_3 = h'(y)t_2$, we find:
\[
B =
\begin{pmatrix}
\frac{C}{\varepsilon_2 \kappa h'(y)} \\
                                          \\
\frac{\varepsilon C}{\varepsilon_2 \kappa} \\
                                            \\
0
\end{pmatrix},
\quad \text{that is}
\]

\begin{align} \label{35}
B_1 = \frac{C}{\varepsilon_2 \kappa h'(y)}, \quad B_2 = \frac{\varepsilon C}{\varepsilon_2 \kappa} \quad \text{ a constant and} \quad B_3 = 0 
\end{align}

By \eqref{24} and \eqref{35}, we have:

\begin{align} \label{36}
\left\lbrace
\begin{array} {llllll}
\langle \nabla_T \partial_y,\; B \rangle = 0 \\
\langle \nabla_T \partial_z,\; B \rangle = -\frac{1}{2} t_3 f_y \left( \frac{\varepsilon C}{\varepsilon_2 \kappa} \right)
\end{array}.
\right.
\end{align}

On the other hand,

\begin{equation*}
AB =
\begin{pmatrix}
0 \\
\varepsilon B_2 \\
B_1
\end{pmatrix},
\quad \text{and} \quad
\gamma(y) =
\begin{pmatrix}
x_0 \\
y \\
h(y)
\end{pmatrix}.
\end{equation*}

So, using \eqref{35} one sees that
\[
(AB)' =
\begin{pmatrix}
0 \\
0 \\
(B_1)'
\end{pmatrix},
\]

and
\begin{align*}
\big(\gamma \bullet (AB)'\big) &= h(y) (B_1)' \cr
&= h(y) \frac{C}{\varepsilon_2 \kappa} \left( \frac{1}{h'(y)} \right)' \cr
&= \frac{C}{\varepsilon_2 \kappa} h(y) f(y) .
\end{align*}

Finally, we find the equation of surface $S$ in the form \eqref{21}:

\begin{eqnarray} \label{A}
z = \varphi(y), \quad \text{where} \quad \varphi(y) = \frac{-2h(y)}{t_3}
\end{eqnarray}

Since
\begin{eqnarray} \label{B}
h'(y) = \frac{1}{\frac{\varepsilon}{2} \int_{y_0}^{y} f(\tau)d\tau + r_0}, 
\end{eqnarray}

and

\[
t_3(y) = \frac{h'(y)}{\sqrt{|\varepsilon + h'(y)^2 f(y)|}} .
\]
\\
The surface $S$ is totally geodesic if and only if the function $f$ satisfies the following differential equation  
\begin{equation} \label{000}
\varphi''(y) +\frac{\varepsilon}{2}\varphi'^3(y)f_y\big(y, \varphi(y)\big) = 0 .
\end{equation}
With \eqref{A} and \eqref{B}.


\begin{thebibliography}{99}
\bibitem{Camara} Camara, E. B., Ndiaye, A., \emph{Family of surfaces with common geodesic in a Walker 3-manifold}, Palest. J. Math., Accepted .

\bibitem{Carmo} Carmo, M. P. do, \emph{Differential geometry of curves and surfaces},  Prentice-Hall, Inc., Englewood Cliffs, N. J., viii+503 (1976) 190-191.

\bibitem {Erjavec} Erjavec, Z., Inoguchi, J-ichi, \emph{Magnetic trajectories in Walker 3-manifold,} J. Math Anal-Appl, 542 (2025) 128806.

\bibitem{Lopez} L\'opez, R., \emph{Differential Geometry of Curves and Surfaces in Lorentz-Minkowski Space.}
Int. Electron. J. Geom., vol 7 (1) (2014) 44-107.

\bibitem{Niang} Niang, A., Ndiaye, A., and Diallo, A. S., \emph{A Classification of Strict Walker 3-Manifold}, Konuralp J. Math., 9 (1) (2021) 148-153.

\bibitem {Walker} Walker, A. G., \emph{Canonical form for a Riemann space with a parallel field of null planes}; Quart. J. Math., Oxford 1 (2) (1950) 69-79.


\end{thebibliography}
\end{document}